\documentclass[12pt]{amsart}
\usepackage{amsmath,amssymb}
\newtheorem{theorem}{Theorem}
\newtheorem{proposition}[theorem]{Proposition}

\newtheorem{corollary}[theorem]{Corollary}
\def\B{\Bbb B}
\def\C{\Bbb C}
\def\D{\Bbb D}
\def\R{\Bbb R}
\def\O{\mathcal O}
\def\Om{\Omega}
\def\eps{\varepsilon}

\title{Estimates of invariant distances on "convex" domains}

\author{Nikolai Nikolov}

\address{Institute of Mathematics and Informatics\\Bulgarian Academy
of Sciences\\ Acad. G. Bonchev 8, 1113 Sofia, Bulgaria\newline
\indent State university of library studies and information
technologies}\email{nik@math.bas.bg}

\subjclass[2010]{32A25, 32F17, 32F45}

\keywords{Carath\'eodory, Kobayashi and Bergman distances, Bergman
and Szeg\"o kernels, convex, convexifiable and $\C$-convex
domains}

\begin{document}

\begin{thanks}{The author would like to thank P.~Pflug, T. Warszawski
and the referee for their remarks on the paper.}
\end{thanks}

\begin{abstract} Estimates for invariant distances of convexifiable,
$\C$-convexifiable and planar domains are given.
\end{abstract}

\maketitle

\section{Introduction and results}

K. Diederich and T. Ohsawa \cite[p. 182]{DO} asked if $D$ is a
smooth bounded pseudoconvex domain in $\Bbb C^n$, then the
following lower bound for the Bergman distance $b_D$ holds: for
fixed $z$ and $w$ close to $\partial D$, one has that
$$b_D(z,w)\ge-c\log d_D(w),$$
where $d_D(w)=\mbox{dist}(w,\partial D)$ and $c>0$ is a constant
depending only on $D.$ Z. B\l ocki \cite[Theorem 1.3]{Blo}
mentioned this fact for bounded convexifiable domains (not
necessarily smooth).

We shall prove the estimate in the case of bounded $\C$-convex
domains (or, more generally, $\C$-convexifiable). Recall that a
set in $\C^n$ is called $\C$-convex if all its intersections with
complex lines are contractible (cf. \cite[p. 25]{APS}). Note that
a $C^1$-smooth domain is $\C$-convex if and only the complex
tangent hyperplane through any boundary point does not intersect
the domain (cf. \cite[Theorem 2.5.2]{APS}).

Let $D$ be a domain in $\C^n.$ Denote by $c_D$ and $l_D$ the
Carath\'eodory distance and the Lempert function of $D,$
respectively:
$$c_D(z,w)=\sup\{\tanh^{-1}|f(w)|:f\in\O(D,\D)\hbox{ with }f(z)=0\},$$
$$l_D(z,w)=\inf\{\tanh^{-1}|\alpha|:\exists\varphi\in\O(\D,D)
\hbox{ with }\varphi(0)=z,\varphi(\alpha)=w\},$$ where $\D$ is the
unit disc (we refer to \cite{JP} for basic properties of the
objects under consideration). The Kobayashi distance $k_D$ is the
largest pseudodistance not exceeding $l_D.$ We have that
$$c_D\le k_D,\quad c_D\le b_D$$
(if $b_D$ is well-defined). Note also that $k_D=l_D$ for any
planar domain $D$ (cf. \cite[Remark 3.3.8(e)]{JP}). By Lempert's
theorem \cite[Theorem 1]{Lem}, combining  with a result by D.
Jacquet \cite[Theorem 5]{Jac}, $c_D=l_D$ on any $C^2$-smooth
bounded $\C$-convex domain $D$ and hence on any convex domain. On
the other hand, it follows by \cite[Theorem 12]{NPZ} that there
exists a constant $c_n>0,$ depending only on $n,$ such that
\begin{equation}\label{comp}k_D\le 4b_D\le
c_nk_D\end{equation} for any $\C$-convex domain $D$ in $\C^n,$
containing no complex lines (then $b_D$ is well-defined). In other
words, to estimate $b_D,$ it is enough to find lower bounds for
$c_D$ and upper bounds for $l_D.$

Recall that $b_D$ is the integrated form of Bergman metric
$$\beta_D(z;X)=\frac{M_D(z;X)}{\sqrt{K_D(z)}}, \quad z\in D,  X\in\C^n,$$
where
$$M_D(z;X)=\sup\{|f'(z)X|:f\in L_h^2(D)\,||f||_D=1,\;f(z)=0\}$$ and
$$K_D(z)=\sup\{|f(z)|^2:f\in L_h^2(D),\;||f||_D\le1\}$$ is the
Bergman kernel on the diagonal ($K_D(z)>0$ is assumed). So,
$$b_D(z,w)=\inf_\gamma\int_0^1\beta_D(\gamma(t);\gamma'(t)),$$
where the infimum is taken over all smooth curves $\gamma:[0,1]\to
D$ with $\gamma(0)=z$ and $\gamma(1)=w.$

Estimates for invariant distances of strictly pseudoconvex domains
in $\C^n$ and pseudoconvex domains of finite type in $\C^2$ can be
found in \cite{BB} (see also \cite{Aba,JP}) and \cite{Her},
respectively.

Recall now in details two estimates. The proof of \cite[Theorem
5.4]{Blo} (cf. also \cite[Proposition 2.4]{Mer}) implies that if
$D$ is a proper convex domain in $\C^n,$ then
\begin{equation}\label{ca}
c_D(z,w)\ge\frac{1}{2}\log\frac{d_D(z)}{d_D(w)}
\end{equation}
(this proof uses only the existence of an appropriate supporting
(real) hyperplane and the formula for the Poincar\'e distance of
the upper half-plane). On the other hand, by \cite[Theorem
1]{NPT}, for any $C^{1+\eps}$-smooth bounded domain there exists a
constant $c>0$ such that
\begin{equation}\label{le}
l_D(z,w)\le-\frac{1}{2}\log(d_D(z)d_D(w))+c
\end{equation}
(see \cite[Proposition 2.5]{FR} for a stronger estimate for
$k_D$).

The smoothness is essential as an example of a $C^1$-smooth
bounded $\C$-convex planar domain shows (see \cite[Example
2]{NPT}). Moreover, using \cite[p. 146, Theorem 7] {SW}, one may
find a bounded $\C$-convex planar domain for which there is no
similar estimate with any constant instead of $-1/2.$

So, it natural to find an upper bound for $l_D$ in the convex case
and a lower bound for $c_D$ in the $\C$-convex case.

\begin{proposition}\label{lem} Let $D$ be a proper convex domain
in $\C^n.$ Then
$$l_D(z,w)\le\frac{||z-w||}{d(z)-d(w)}\log\frac{d(z)}{d(w)}
\le\frac{||z-w||}{\min(d(z),d(w))}.\footnote{If $d(z)=d(w),$ then
$l_D(z,w)\le||z-w||/d(w).$}$$ In particular, if, in addition, $D$
is bounded, then for any compact subset $K$ of $D$ there is a
constant $c_K>0$ such that
$$b_D\le-c_K\log d_D(w)+1/c_K,\quad z\in K,w\in D.$$
\end{proposition}

The last estimate for $k_D$ instead of $b_D$ (and $K$ a singleton)
is the content of \cite[Proposition 2.3]{Mer}. Similar estimates
for the Kobayashi distance of pseudoconvex Reinhardt domains can
be found in \cite{War}.

\begin{proposition}\label{car} Let $D$ be a proper $\C$-convex domain
in $\C^n.$ Then
$$c_D(z,w)\ge\frac{1}{4}\log\frac{d_D(z)}{4d_D(w)}.$$
Hence, if, in addition, $D$ is bounded, then for any compact
subset $K$ of $D$ there is a constant $c_K>0$ such that
$$b_D(z,w)\ge-\frac{1}{4}\log d_D(w)-c_K,\quad z\in K,w\in D.$$
\end{proposition}

Note that by \cite[p. 2381]{CF} the first estimate in Proposition
\ref{car} implies the following

\begin{corollary}\label{sze}
The Bergman and Szeg\"o kernels (on the diagonal) are comparable
on any $C^2$-smooth bounded $\C$-convex domain.
\end{corollary}

We point out that \cite[Theorem 1.3]{CF} deals with the convex
case.
\smallskip

\noindent\textsc{Remark.} (a) The estimate for $l_D$ is sharp when
$z\to w.$ Moreover, it is sharp up to a constant when $z$ is fixed
and $w\to\partial D.$ Indeed, denote by $R_D(z,w)$ the right-hand
side of the first inequality in Proposition 1. If
$\theta\in(0,\pi)$ and $D_\theta=\{z\in\C_\ast:|\arg z|<\theta\},$
then
$$\lim_{\theta\to 0}\lim_{x\to
0+}\frac{l_{D_\theta}(1,x)}{R_{D_\theta}(1,x)}=\frac{\pi}{4}.$$

\noindent (b) The factor $1/4$ in the bound for $c_D$ is optimal
as $D=\C_\ast\setminus\Bbb R^+$ shows.
\smallskip

\noindent (c) Estimates for the infinitesimal forms of the
distances under consideration, namely, the Carath\'eodory,
Kobayashi and Bergman metrics, of convex and $\C$-convex domains
can be found in \cite{NPZ}. The bounds there depend only on the
distance to the boundary from the respective point in the
respective direction.
\smallskip

Our main result is in the spirit of \cite[Theorem 1.3]{Blo}, where
a lower bound for the Bergman metric is mentioned in the locally
convexifiable case (and a hint for a proof is given).

\begin{proposition}\label{kob} Let $D$ be a bounded domain in
$\C^n$ which is locally $\C$-convexifiable, i.e. for any point
$a\in\partial D$ there exist a neighborhood $U_a$ of $a,$ an open
set $V_a$ in $\C^n$ and a biholomorphism $F_a:U_a\to V_a$ such
that $F_a(D\cap U_a)$ is $\C$-convex. Then there exists a constant
$c>0$ such that for any compact subset $K$ of $D$ one can find a
constant $c_K>0$ with
$$s_D(z,w)\ge-c\log d_D(w)-c_K,\quad z\in
K,w\in D,$$ where $s_D=k_D$ or $s_D=b_D.$

Moreover, if $D$ is locally convexifiable or $C^{1+\eps}$-smooth
and locally $\C$-confexifiable, then for any compact subset $K$ of
$D$ one can find a constant $c'_K>0$ with
$$s_D(z,w)\le-c_K'\log d_D(w)+1/c'_K,\quad
z\in K,w\in D.$$
\end{proposition}

Finally, we consider the planar case. We shall say that  a
boundary point $p$ of a planar domain $D$ is Dini-smooth if
$\partial D$ near $p$ is a Dini-smooth curve
$\gamma:[0,1]\to\C$.\footnote {This means that
$\int_0^1\frac{\omega(t)}{t}dt<\infty,$ where $\omega$ is the
modulus of continuity of $\gamma'.$} We shall that a planar domain
is Dini-smooth if it is Dini-smooth near an boundary point.

\begin{proposition}\label{pla} Let $p$ be a Dini-smooth
boundary point of a planar domain $D.$ Then for any neighborhood
$U$ of $p$ and any compact subset $K$ of $D$ there exist a
neighborhood $V$ of $p$ and a constant $c>0$ such that
$$s_D(z,w)\ge-\frac{1}{2}\log d_D(w)-c,\quad z\in D\setminus U, w\in D\cap V,$$
$$|s_D(z,w)+\frac{1}{2}\log d_D(w)|\le c,\quad z\in K, w\in D\cap V,$$
where $s_D=c_D,$ $s_D=l_D(=k_D)$ or $s_D=b_D/\sqrt 2.$
\end{proposition}

Since $k_D$ and $b_D$ are the integrated forms of $\kappa_D$ and
$\beta_D,$ we get the following

\begin{corollary}\label{cor} Let $p$ and $q$ be different Dini-smooth
boundary points of a planar domain $D.$ If $s_D=l_D(=k_D)$ or
$s_D=b_D/\sqrt 2,$ then the function $$2s_D(z,w)+\log d_D(z)+\log
d_D(w)$$ is bounded for $z$ near $q$ and $w$ near $p.$
\end{corollary}

In general, $c_D$ is not an inner distance (even in the plane).
So, the next proposition is not a direct consequence of
Proposition \ref{pla}.

\begin{proposition}\label{cara} Let $p$ and $q$ be different Dini-smooth
boundary points of a planar domain $D.$ Then the function
$$2c_D(z,w)+\log d_D(z)+\log d_D(w)$$ is bounded for $z$ near $q$
and $w$ near $p.$
\end{proposition}

The next result is optimal for the boundary behavior of $c_D$ and
$l_D(=k_D)$ in the planar case. It is more general than the last
results but its proof use these results. Similar (and slightly
weaker) result for $k_D$ on $C^2$-smooth strictly pseudoconvex
bounded follows by \cite[Theorem 1, Proposition 1.2]{BB}.

\begin{proposition}\label{gen} Let $D$ be a Dini-smooth bounded planar
domain.\footnote{This means that $D$ is Dini-smooth near any
boundary point.} Then there exists a constant $c\ge 1$ such that
$$\log\left(1+\frac{|z-w|}{c\sqrt{d_D(z)d_D(w)}}+\frac{|z-w|^2}{cd_D(z)d_D(w)}
\right)\le 2c_D(z,w)$$ $$\le 2l_D(z,w)\le
\log\left(1+\frac{c|z-w|}{\sqrt{d_D(z)d_D(w)}}+\frac{c|z-w|^2}{d_D(z)d_D(w)}\right).$$

In particular, the function $l_D-c_D$ is bounded on $D\times D.$
\end{proposition}

It is shown in \cite[Theorem 1]{Ven} that if $D$ is strongly
pseudoconvex domain in $\C^n,$ then
$$\lim_{\substack{w\to\partial D\\z\neq
w}}\frac{c_D(z,w)}{k_D(z,w)}=1\quad\mbox{uniformly in }z\in D.$$
We have the following planar extension of this result.

\begin{proposition}\label{ven} If $D$ is finitely connected bounded planar domain
without isolated boundary points, then
$$\lim_{\substack{w\to\partial D\\z\neq
w}}\frac{c_D(z,w)}{l_D(z,w)}=1\quad\mbox{uniformly in }z\in D.$$
\end{proposition}

\section{Proofs}

\noindent{\it Proof of Proposition \ref{lem}.} Denote by $C_{z,w}$
the convex hull of the union of the discs $\D(z,d_D(z))$ and
$\D(w,d_D(w)),$ lying in the complex line through $z$ and $w.$ Let
$\gamma(t)=z+t(w-z).$ Since $C_{z,w}\subset D$ and
$l_{C_{z,w}}=k_{C_{z,w}}$ is the integrated form of the Kobayashi
metric $\kappa_{C_{z,w}},$ \footnote{If $D\subset\C^n,$ then
$\kappa_D(z;X)=\inf\{|\alpha|:\exists\varphi\in\O(\D,D) \hbox{
with }\varphi(0)=z,\alpha\varphi'(0)=X\}.$} then
$$l_D(z,w)\le
l_{C_{z,w}}(z,w)\le\int_0^1\kappa_{C_{z,w}}(\gamma(t);\gamma'(t))dt$$
$$\le\int_0^1\frac{|\gamma'(t)|}{d_{C_{z,w}}(\gamma(t))}dt=
\frac{||z-w||}{d(z)-d(w)}\log\frac{d(z)}{d(w)}.$$

This inequality and \eqref{comp} lead to the wanted result for
$b_D.$
\smallskip

\noindent{\it Proof of Proposition \ref{car}.} Let
$p(w)\in\partial D$ be such that $||w-p(w)||=d_D(w).$ Since $E$ is $\Bbb C$-convex,
there exists a hyperplane $H_{p(w)}$ through $p(w)$ and disjoint from $D$ (cf. \cite
[Theorem 2.3.9(ii)]{APS}.). Denote by
$D_w$ and $z_w$ be the projections of $D$ and $z$ onto
the complex line through $w$ and $p(w)$ in direction $H_{(p(w)},$ respectively. By
\cite[Theorem 2.3.6]{APS}, $D_w$ is a simply connected domain and $p(w)\in\partial D_w.$
Denote by $\psi_w\in\O(\D,D_w)$ a Riemann map such that
$\psi_w(0)=z_w.$ If $\psi_w(\alpha_w)=w,$ then
$$c_D(z,w)\ge c_{D_w}(z_w,w)=\tanh^{-1}|\alpha_w|.$$ By \cite[p.
139, Corollary 6]{SW} (which is a consequence of the K\"obe 1/4
and the K\"obe distortion theorems),
$$\tanh^{-1}|\alpha_w|\ge\frac{1}{4}\log\frac{|\psi_w'(0)|}{4d_{D_w}(w)}.$$
Since $d_{D_w}(w)=d_D(w)$ and $|\psi_w'(0)|\ge d_{D_w}(z_w)\ge
d_D(z),$ it follows that
$$c_D(z,w)\ge\frac{1}{4}\log\frac{d_D(z)}{4d_D(w)}.$$

This inequality and $b_D\ge c_D$ imply the desired result for
$b_D.$
\smallskip

\noindent{\it Proof of Proposition \ref{kob}.}\footnote{Some
difficulty arises from the fact that, in contrast to invariant
metrics, general localization principles for invariant distances
are not known. However, a strong localization principle holds for
$k_D$ and $c_D$ if $D$ is strongly pseudoconvex (see
\cite[Proposition 3, Theorem 1]{Ven}.} First, we shall prove the
lower bound.

Note that \begin{equation}\label{est} 0<c_a\le\frac{d_{F_a(D\cap
U_a)}(F_a(w))}{d_D(w)}\le\frac{1}{c_a}\quad\mbox{near any
}a\in\partial D. \end{equation} Then, by Proposition \ref{car}, we
may find a finite set $M\subset\partial D$ and a constant $c_1>0$
such that
$$s_{D\cap U_a}(z,w)\ge\frac{1}{4}\log\frac{d_D(z)}{d_D(w)}-c_1,\quad
z,w\in D\cap V_a, a\in M,$$ where $V_a\subset U_a$ is a
neighborhood of $a$ such that $\partial D\subset\cup_{a\in M}
V_a.$

Denote now by $S_D$ the Kobayashi or Bergman metrics of $D.$ By
localization principles (cf. \cite[Proposition 7.2.9 and
Proposition 6.3.5]{JP}, since $D$ is pseudoconvex), there exists a
constant $c_2>0$ such that
$$S_D\ge 4c_2 S_{D\cap U_a}\hbox{ on }(D\cap V_a)\times\C^n.$$

Let $W_a\Subset V_a$ be such that $W=\cup_{a\in M} W_a$ does not
intersect $K$ and contains $\partial D.$  Set $r=\min_{a\in
M}\mbox{dist}(\partial W_a,\partial V_a).$

Let $\eps>0.$ Since $s_D$ is the integrated form of $S_D,$ for any
$z\in K$ and $w\in D\cap W$ there exists a smooth curve
$\gamma:[0,1]\to D$ with $\gamma(0)=z,$ $\gamma(1)=w$ and
$$s_D(z,w)+\eps>\int_0^1S_D(\gamma(t);\gamma'(t))dt.$$
Let $t_1=\max\{t\in(0,1):\gamma(t)\in G=D\setminus W\}.$ Choose a
point $a_1\in M$ such that $\B_n(\gamma(t_1),r)\subset V_{a_1}.$
Let $t_2=\sup\{t\in(t_1,1]:\gamma([t_1,t))\in V_{a_1}\}$ and etc.
In this way we may numbers $0<t_1<\dots<t_{N+1}=1$ and points
$a_1,\dots,a_{N+1}\in M$ such that $\gamma[t_j,t_{j+1})\subset
D\cap V_{a_j}$ and $||\gamma(t_{j+1})-\gamma(t_j)||\ge r,$ $1\le
j\le N.$ Then
$$s_D(z,w)+\eps>c_2\sum_{j=1}^N s_{D\cap
U_{a_j}}(\gamma(t_j),\gamma(t_{j+1}))$$
$$\ge
c_2\sum_{j=1}^N\log\frac{d_D(\gamma(t_j))}{d_D(\gamma(t_{j+1}))}
-c_3N$$
$$\ge c_2\log\frac{\mbox{dist}(G,\partial D)}{d_D(w)}-c_3N,$$
where $c_3=4c_1c_2.$

On the other hand, since $D$ is a bounded domain, there exists a
constant $c_4>0$ such that $s_D(z_1,z_2)\ge c_4||z_1-z_2||.$ Then
$$s_D(z,w)+\eps>\sum_{j=1}^Ns_D(\gamma(t_j),\gamma(t_{j+1}))\ge c_4rN.$$

So,
$$\left(1+\frac{c_3}{c_4r}\right)(s_D(z,w)+\eps)\ge c_2\log\frac{\mbox{dist}
(G,\partial D)}{d_D(w)}.$$

The case when $w\in G$ is trivial which completes the proof of the
lower bound.

The proof of the upper bound is easier. Fix a point $a\in\partial
D.$ It is enough to find a constant $c'_{a,K}>0$ such that the
estimate holds for $w$ near $a.$ Take a point $u\in U_a$ and a
neighborhood $V_a\Subset U_a$ of $a$ and a point $u\in D\cap U_a.$
It follows by Proposition \ref{lem}, \eqref{le} and \eqref{est}
that
$$k_D(z,w)\le k_D(z,u)+k_D(u,w)\le k_D(z,u)+k_{D\cap U}(u,w)$$
$$\le 1/c'_{a,K}-c'_{a,K}\log d_D(w),\ z\in K,w\in D\cap V_a.$$

The upper bound for $b_D$ follows similarly. It suffices to use
that
$$b_D\le\widetilde{c}_ab_{D\cap U_a}\le\widetilde{c}_ac_n k_{D\cap
U_a}$$ in view of \cite[Proposition 6.3.5]{JP} and \eqref{comp}.
\smallskip

\noindent{\it Proof of Proposition \ref{pla} for $c_D$ and $l_D$.}
We may find a Dini-smooth Jordan curve $\zeta$ such that
$\zeta=\partial D$ near $p$ and $D\subset G:=\zeta_{\mbox{ext}}.$
Take a point $a\not\in\overline{G}$ and consider the union $G_e$
of $0$ and the image of $G$ under the map
$\varphi:z\to(z-a)^{-1}.$ There exists a conformal map
$\psi:G_e\to\D.$ It extends to a $C^1$-diffeomorphism from
$\overline{G_e}$ to $\overline{\D}$ (cf. \cite[Theorems
3.5]{Pom}). Setting $\eta=\psi\circ\varphi,$ then
$$c_D(z,w)\ge c_\D(\eta(z),\eta(w)).$$
Now the lower bound for $c_D$ follows by the same bound for $c_\D$
and an inequality of type \eqref{est}.

The estimate
$$l_D(z,w)\le-\frac{1}{2}\log d_D(w)-c,\quad z\in K, w\in D\cap V$$
follows by \eqref{le}. It can be also obtained in the following
way. There exist a Dini-smooth domain simply connected domain
$G_i\subset D$ and a neighborhood $V$ of $p$ such that $\partial
G\cap V=\partial D\cap V.$ Take a point $u\in V.$ Since $l_D=k_D,$
then $$k_D(z,w)\le k_D(z,u)+k_{G_i}(u,w).$$ It remains to repeat
the final arguments from the first paragraph.
\smallskip

\noindent{\it Proof of Proposition \ref{pla} for
$b_D.$}\footnote{We have to modify the previous proof, since the
Bergman distance is not monotone under inclusion of planar
domains; to see this, use \cite[Example 7]{PZ}.} Choosing $G$ as
above, then
$$b_D(z,w)=b_{\eta(D)}(\eta(z),\eta(w)).$$ By the
Dini-smoothness,
$$\lim_{z\to p}\frac{d_{\eta(D)}(\eta(z))}{d_D(z)}=|\eta'(p)|.$$
We may assume that $\eta(p)=1.$ So, it is enough to get the
estimates for $D\subset\D$ such that $F=\D\cap\D(1,r)\subset D$
for some $r\in(0,1).$

First, we shall prove that if $0<r'<r,$ then
$$\sqrt 2b_D(z,w)\le-\log d_D(w)+c',\quad z\in K, w\in F'=\D\cap\D(1,r')$$
for some constant $c'>0.$

For a domain $\Om\subset\C$ set $\beta_\Om(z)=B_\Om(z;1)$ and
$\kappa_\Om(z)=\kappa_\Om(z;1).$ Let $\check F=\D\setminus F$ and
$$l_\D(u,\check F)=\inf_{w\in\check F}l_\D(u,w).$$

Then for any $r''\in(r',r)$ we may find a constant $\tilde c>0$
such that
$$\beta_D(u)\le\beta_F(u)\sqrt\frac{K_F(u)}{K_\D(u)}
=\frac{\sqrt2\kappa^2_F(u)}{\kappa_\D(u)}$$
$$\le\sqrt 2\coth^2l_\D(u,\check F)\kappa_\D(u)\le\frac{\sqrt 2}{1-|u|^2}+\tilde c,
\quad u\in F''=\D\cap\D(1,r'').$$ (for the equality use that $F$ is biholomorphic to
$\D$ and for the inequality "between the lines" cf. \cite[Proposition 7.2.9]{JP}).

Let $z\in K,$ $w\in F'$ and $w'=[0,w]\cap\partial D(1,r'').$ Then
$$b_D(z,w)\le b_D(z,w')+|w-w'|\left(\tilde c+\sqrt
2\int_0^1\frac{dt}{1-|w'+t(w-w')|^2}\right)$$
$$\le(-\log d_D(w)+c')/\sqrt2$$ for some constant $c'>0.$

Now, shrinking $r$ such that $\D(1,r)\subset U,$ it remains to
prove that
$$\sqrt 2b_D(z,w)\ge-\log d_D(w)-c'',\quad z\in\check F, w\in F'$$
for some constant $c''>0.$

We have that
$$\beta_D(u)\ge\beta_\D(u)\sqrt\frac{K_\D(u)}{K_F(u)}
=\frac{\sqrt2\kappa^2_\D(u)}{\kappa_F(u)}$$
$$\ge\sqrt 2\tanh l_\D(u,\check F)\kappa_\D(u)\ge\frac{\sqrt 2}{1-|u|^2}-\hat c,
\quad u\in F''.$$

For $z\in\check F,$ $w\in F'$ and $\eps>0$ there exists a smooth
curve $\gamma:[0,1]\to D$ with
$$b_D(z,w)+\eps>\int_0^1\beta_D(\gamma(t))|\gamma'(t)|dt.$$
Let $t_0=\sup\{t\in(0,1):\gamma(t)\not\in F''\}$ and
$$\beta^{\hat c}_\D(z;X)=|X|\left(\frac{\sqrt 2}{1-|z|^2}-\hat c\right)^+.$$
Then
$$b_D(z,w)+\eps>\int_{t_0}^1b_D(\gamma(t))|\gamma'(t)|dt\ge
\hat b_\D(w,\check F),$$ where $\hat b$ is the integrated form of
the Finsler pseudometric
$$\hat\beta_\D(u;X)=|X|\left(\frac{\sqrt 2}{1-|u|^2}-\hat c\right)^+.$$
It remains to use that, shrinking $r'$ (if necessary),
$$\hat b_\D(w,\check F)\ge(-\log d_D(w)-c'')/\sqrt2$$ for some constant $c''>0$
(cf. \cite[Theorem 1.1]{BB}).
\smallskip

\noindent{\it Proof of Corollary \ref{cor}}. Since $k_D$ and $b_D$
are the integrated forms of $\kappa_D$ and $b_D,$ the boundedness
from below follows by the first inequality in Proposition
\ref{pla} (cf. the proof of \cite[Proposition 10.2.6]{JP}).
Choosing a point $a\in D,$ the boundedness from above is a
consequence of the inequality $s_D(z,w)\le s_D(z,a)+s_D(a,w)$ and
the second inequality in Proposition \ref{pla}.
\smallskip

\noindent{\it Proof of Proposition \ref{cara}}. In virtue of the
inequality $c_D\le k_D$ and Corollary \ref{cor}, we have to prove
only the boundedness from below. For this, take disjoint
Dini-smooth Jordan curves $\zeta'$ and $\zeta''$ such that
$\zeta'=\partial D$ near $p,$ $\zeta''=\partial D$ near $q$ and
$D\subset G:=\zeta'_{\mbox{ext}}\cap\zeta''_{\mbox{ext}}.$ Note
that any Dini-smooth bounded double connected planar $\tilde{G}$
domain can be conformally map to some annulus
$A_r=\{z\in\C:1/r<|z|<r\}$ ($r>1$) and the respective mapping
extends to a $C^1$-diffeomorphism from $\overline{\tilde G}$ to
$\overline{A_r}.$\footnote{To see this, we can proceed as follows
(S. R. Bell, private communication). First, take a conformal
mapping $\varphi_1$ from the domain bounded by the outher boundary
of $\tilde G$ to $\D.$ Next, choose a point $a$ in the interior of
the inner boundary $\Gamma$ of $\psi_1(\tilde G)$ and set
$\psi_2:z\to(z-a)^{-1}.$ Consider now a conformal mapping $\psi_3$
from the domain bounded by $\psi_2(\Gamma)$ to $\D.$ Then
$\psi=\psi_3\circ\psi_2\circ\psi_1$ maps conformally $\tilde G$ to
a bounded double connected planar domain $G'$ with real-analytic
boundary. It remains to apply the reflection principle to a
conformal mapping from $G'$ to $A_r.$}

Then, proceeding similarly to the proof of Proposition \ref{pla}
for $c_D,$ it is enough to show that
$$2c_{A_r}(z,w)+\log d_{A_r}(z)+\log d_{A_r}(w)$$ is bounded from below
for $z\in\R$ near $r$ and $w$ near $p,$ where $|p|=1/r;$ this is
equivalent to
$$m_{A_r}(z,w):=\tanh c_{A_r}(z,w)\ge 1-c d_{A_r}(z)d_{A_r}(w)$$
for some constant $c>0.$

Recall that (cf. \cite[Proposition 5.5]{JP})
$$m_{A_r}(z,w)=\frac{f(z,w)f(1/z,-|w|)}{r|w|},$$
where $f$ is a holomorphic function on $\overline{A_r\times
A_r}\setminus\{u=v\in\partial A_r\}$ and $|f(u,v)|=1$ if $|u|=r,$
$v\in\overline{A_r}$ or $u\in\overline{A_r},$ $|v|=1/r$ ($u\neq
v$).

In particular, $$\frac{\partial^n f}{\partial
u^n}=\frac{\partial^n f}{\partial v^n}=0,\quad n\in\Bbb N,$$ at
any point $(u,v)$ with $|u|=r$ and $|v|=1/r.$ Then, by the Taylor
expansion,
$$|f(z,w)-f(r,w(r|w|)^{-1})|\le c_1d_{A_r}(z)d_{A_r}(w).$$
This implies that
$$|f(z,|w|)-f(r,1/r)|\le c_1d_{A_r}(z)d_{A_r}(|w|)$$
(the constant can be chosen the same for $z$ near $r$ and $w$ away
from $r$). Since $f(r,\cdot)$ is a unimodular constant and
$d_{A_r}(w)=d_{A_r}(|w|),$ it follows that
$$|m_{A_r}(z,w)-m_{A_r}(z,|w|)|\le c_2d_{A_r}(z)d_{A_r}(|w|).$$

Further, $c_{A_r}(z,|w|)=c_{A_r}(z,t)+c_{A_r}(t,|w|)$ for
$t\in[|w|,z]$ (cf. \cite[Lemma 5.11(b)]{JP}). Then Proposition
\ref{pla} implies that
$$m_{A_r}(z,|w|)\ge 1-c_3d_{A_r}(z)d_{A_r}(|w|).$$
Hence we may choose $c=c_2+c_3$ which completes the proof.
\smallskip

\noindent{\it Proof of Proposition \ref{gen}}. Using Corollary
\ref{cor} and Proposition \ref{cara}, it is enough to prove the
inequalities for $z$ and $w$ near a fixed point $p\in\partial D.$
Moreover, it is easy to see that these inequalities are equivalent
to $$\frac{|z-w|}{\sqrt{cd_D(z)d_D(w)+|z-w|^2}}\le\tanh c_D(z,w)$$
$$\le\tanh
l_D(z,w)\le\frac{|z-w|}{\sqrt{c^{-1}d_D(z)d_D(w)+|z-w|^2}}$$ for
some constant $c\ge 1$.\footnote{These estimates implies the
bounds for the Green function $g_D$ from the crucial Lemma 4.2 in
\cite{Swe}, since $\tanh c_D\le\exp(-2\pi g_D)\le \tanh l_D.$}

To prove the lower bound for $\tanh c_D(z,w),$ let $\eta$ be as in
the proof of Proposition \ref{pla} for $c_D$ and $l_D$. Then it is
not difficult to find a constant $c_1>0$ such that
$$\tanh c_D(z,w)\ge\tanh c_\D(z_1,w_1)\ge
\frac{|z_1-w_1|}{\sqrt{c_1d_\D(z_1)d_\D(w_1)+|z_1-w_1|^2}},$$
where $z_1=\eta(z)$ and $w_1=\eta(w).$ It remains to use that,
similarly to \eqref{est}, $d_D\ge c_2 d_\D$ and $|z_1-w_1|\ge
c_2|z-w|$ for some constant $c_2>0.$

The proof of the upper bound for $\tanh l_D(z,w)$ is similar (by
using $G_i$ from the second part of the proof mentioned above) and
we skip it.
\smallskip

\noindent{\it Proof of Proposition \ref{ven}}. By the K\"obe
uniformization theorem, we may assume that $\partial D$ consists
of disjoint circles. Using Proposition \ref{gen} and compactness,
it is enough to prove that for any point $p\in\partial D,$
$$\lim_{z\neq w\to p}\frac{c_D(z,w)}{l_D(z,w)}=1.$$
Applying an inversion, we may suppose that the outer boundary of
$D$ is the unit circle $\Gamma$ and $p\in\Gamma.$ Let $U$ be a
disc centered at $p$ such that $\D\cap U\subset D.$ Then
$$1\ge \frac{c_D(z,w)}{l_D(z,w)}\ge\frac{c_\D(z,w)}
{l_{\D\cap U}(z,w)}=\frac{k_\D(z,w)}{k_{\D\cap U}(z,w)}.$$
Considering $\D$ as a part of the unit ball in $\C^2,$ it follows
that the last ratio tends to 1 as a particular case of the same
result for strongly pseudoconvex domains (see \cite[Proposition
3]{Ven}).


\begin{thebibliography}{}

\bibitem{Aba} M. Abate, {\it  Iteration theory of holomorphic maps on taut
manifolds,} Mediterranean Press, 1989.

\bibitem{APS} M. Andersson, M. Passare, R. Sigurdsson, {\it Complex
convexity and analytic functionals,} Birkh\"auser, 2004.

\bibitem{BB} Z. M. Balogh, M. Bonk, {\it Gromov hyperbolicity and
the Kobayashi metric on strictly pseudoconvex domains}, Comment.
Math. Helv. 75 (2000) 504-533.

\bibitem{Blo} Z. B\l ocki, {\it The Bergman metric and the
pluricomplex Green function}, Trans. Amer. Math. Soc. 357 (2005),
2613-2625.

\bibitem{CF} B.-Y. Chen, S. Fu, {\it Comparison of the Bergman and
Szeg\"o kernels}, Adv. Math. 228 (2011), 2366-2384.

\bibitem{DO} K. Diederich, T. Ohsawa, {\it An estimate for the Bergman
distance on pseudoconvex domains}, Ann. Math. 141 (1995), 181-190.

\bibitem{FR} F. Forstneri\v c, J. P. Rosay,
{\it Localization of the Kobayshi metric and the boundary
continuity of proper holomorphic mappings}, Math. Ann 279 (1987),
239-252.

\bibitem{Her} G. Herbort, {\it Estimation on invariant distances on
pseudoconvex domains of finite type in dimension two}, Math. Z.
251 (2005), 673-703.

\bibitem{Jac} D. Jacquet, {\it $\Bbb C$-convex domains with $C^2$ boundary},
Complex Var. Elliptic Equ. 51 (2006), 303-312.

\bibitem{JP} M. Jarnicki, P. Pflug, {\it Invariant distances and metrics
in complex analysis}, de Gruyter Exp. Math. 9, de Gruyter, Berlin
(1993).

\bibitem{Lem} L. Lempert, {\it Intrinsic distances and holomorphic retracts,}
Complex analysis and applications '81, Sofia (1984), 341-364.

\bibitem{Mer} P. R. Mercer, {\it Complex geodesics and iterates of holomorphic
maps on convex domains in $\C^n$}, Trans. Amer. Math. Soc. 338
(1993), 201-211.

\bibitem{NPT} N. Nikolov, P. Pflug, P. J. Thomas, {\it Upper bound
for the Lempert function of smooth domains}, Math. Z. 266 (2010),
425-430.

\bibitem{NPZ} N. Nikolov, P. Pflug, W. Zwonek, {\it Estimates for
invariant metrics on $\Bbb C$-convex domains}, Trans. Amer. Math.
Soc. 363 (2011), 6245-6256.

\bibitem{PZ} P. Pflug, W. Zwonek, {\it Logarithmic capacity and Bergman
functions}, Arch. Math. 80 (2003), 536-552.

\bibitem{SW} W. Seidel, J. L. Walsh, {\it On the derivatives of functions
analytic in the unit circle and their radii of univalence and of
p-valence}, Trans. Amer. Math. Soc. 52 (1942), 128-216.

\bibitem{Swe} G. Sweers, {\it Positivity for a Strongly Coupled Elliptic
System by Green function Estimates}, J. Geom. Anal. 4 (1994),
121-142.

\bibitem{Ven} S. Venturini, {\it Comparision between the Kobayashi
and Carath\'eodory distances on strongly pseudoconvex bounded
domains in $\C^n$}, Proc. Amer. Math. Soc. 107 (1989), 725-730.

\bibitem{War} T. Warszawski, {\it Boundary behavior of the Kobayashi distance
in pseudoconvex Reinhardt domains}, Michigan Math. J. 61 (2012),
575-592.

\bibitem{Pom} Ch.~Pommerenke, {\it Boundary behaviour of
conformal maps}, Grundl. math. Wissensch. 299, Springer, Berlin
(1992).

\end{thebibliography}
\end{document}